\documentclass[aap,MSNbibl,seceqn,dvips]{arximspdf}
\usepackage{mathbh}

\doi{10.1214/13-AAP972} 
\volume{24}
\issue{5}
\pubyear{2014}
\firstpage{2070}
\lastpage{2090}

\makeatletter
\newcommand{\eqref}[1]{(\ref{#1})}
\newtheorem{theorem}{Theorem}[section]
\newtheorem{lemma}[theorem]{Lemma}

\newtheorem{proposition}[theorem]{Proposition}

\newproclaim{definition}[theorem]{Definition}
\newproclaim{assumption}[theorem]{Assumption}

\newproclaim{remark}[theorem]{Remark}

\def\N{\mathbb{N}}
\def\R{\mathbb{R}}

\def\ep{\varepsilon}

\def\E{\mathbf{E}}
\def\EE{\mathbb{E}}
\def\P{\mathbf{P}}
\def\PP{\mathbb{P}}
\def\F{\mathcal{F}}

\def\Ind{\mathbh{1}}

\def\to{\rightarrow}

\newcommand{\dd}{\mathrm{d}}

\def\root{\rho}
\def\Lt{\mathrm L}
\def\Lc{\mathcal L}
\def\widebar{\overline} 

\newcommand{\parent}{\vec}

\makeatother

\begin{document}
\begin{frontmatter}

\title{Performance of the Metropolis algorithm on a~disordered tree: The
Einstein relation}
\runtitle{Metropolis algorithm on a disordered tree}

\begin{aug}
\author[A]{\fnms{Pascal} \snm{Maillard}\corref{}\thanksref{t1}\ead[label=e1]{pascal.maillard@weizmann.ac.il}}
\and
\author[A]{\fnms{Ofer} \snm{Zeitouni}\thanksref{t2}}
\thankstext{t1}{Supported in part by a grant from the Israel Science
Foundation.}
\thankstext{t2}{Supported in part by NSF Grant DMS-12-03201,
the Israel Science Foundation and the Herman~P. Taubman chair
of Mathematics at the Weizmann Institute.}
\runauthor{P. Maillard and O. Zeitouni}
\affiliation{Weizmann Institute of Science}
\address[A]{Department of Mathematics and Computer Science\\
Weizmann Institute of Science\\
P.O. Box 26\\
Rehovot 76100\\
Israel\\
\printead{e1}} 
\end{aug}

\received{\smonth{4} \syear{2013}}
\revised{\smonth{7} \syear{2013}}

%
\begin{abstract}
Consider a $d$-ary rooted tree ($d\geq3$) where each edge $e$
is assigned
an i.i.d. (bounded) random variable $X(e)$ of negative mean. Assign to
each vertex $v$ the sum $S(v)$ of $X(e)$ over all edges connecting $v$
to the
root, and assume that the maximum $S_n^*$ of $S(v)$ over all vertices
$v$ at
distance $n$ from the root tends to infinity (necessarily, linearly) as
$n$ tends to infinity.
We analyze the Metropolis algorithm on the tree and show that under
these assumptions there always exists a temperature $1/\beta$ of the algorithm
so that it achieves a linear (positive) growth rate in linear time.
This confirms a conjecture of Aldous
[\textit{Algorithmica} \textbf{22} (1998) 388--412].
The proof is obtained
by establishing an Einstein relation for the Metropolis algorithm on
the tree.
\end{abstract}

%

\end{frontmatter}

\section{Introduction}
Given a $d$-regular rooted tree, attach to each edge $e$ a random
variable $X(e)$, such that the variables are independent and
identically distributed. For a vertex $v$ in the tree, denote by $S(v)$
the sum of the variables $X(e)$ over all edges $e$ on the path from the
root to $v$. This defines a branching random walk, a basic model for a
disordered tree. It is natural to ask for an efficient algorithm which
explores the vertices of this tree in order to find vertices $v$ with a
large value of $S(v)$. In fact, Aldous \cite{Aldous1998} proposed this
problem as a benchmark problem for comparing different generic
optimization algorithms, since the na\"{\i}ve
approach, which would be to simply explore all
vertices down to the level $n$
in the tree and taking the one with the maximal value of $S(v)$, is
a bad choice for an algorithm because
the number of vertices grows exponentially in $n$.

The Metropolis algorithm is a general recipe for constructing a
discrete-time Markov chain on a finite state space for which a given
distribution $\pi$ is stationary and whose transitions respect a given
graph structure of the state space. In the context of the comparison of
algorithms
discussed earlier, Aldous \cite{Aldous1998} suggested using the
Metropolis algorithm to ``sample'' a certain Gibbs measure on the
vertices of a branching random walk tree, namely the one which assigns
mass $e^{\beta S(v)}$ to a vertex~$v$, for some parameter $\beta>0$.
In the case where this measure is infinite, for example when there is
an infinite number of vertices $v$ with $S(v)\ge0$, this algorithm
should ``walk down the tree'' and, for an appropriate choice of the
parameter~$\beta$, find vertices $v$ with high values of $S(v)$. Let
$|v|$ denote the level of the vertex $v$ in the tree, and let $V_k$ be
the vertex visited by the Metropolis algorithm at the time~$k$. Aldous
raised the following natural
question: If the maximum of the branching random walk has positive
speed, that is, if $\lim_{n\to\infty} \max_{|v|=n} S(v)/n > 0$, does
there exist a choice of the parameter $\beta$, such that $\liminf
_{k\to
\infty} S(V_k)/k > 0$? We will answer this question in the affirmative
for a certain class of laws of the variables $X(e)$, including the
binomial distribution.

In fact, we show more: Let $v_\beta= \lim_{k\to\infty} S(V_k)/k$,
which exists almost surely~\cite{Aldous1998}.
We show that there exists a parameter $\beta_0>0$, such that $v_{\beta_0}=0$
and $(\dd v_\beta/\dd\beta) |_{\beta=\beta_0} = \sigma^2/2$, where
$\sigma^2$ is the asymptotic variance of $S(V_k)$, which we show to be
positive and finite. This result was conjectured by Aldous, who gave
heuristic arguments and numerical evidence for it [in the case where
the variables $X(e)$ only take the values $1$ and $-1$]. Results of
this type are also known as \emph{Einstein relations} in the domain of
random walks in random environments and our methods of proof will
indeed rely on many techniques from this field, some of which have been
obtained recently.

\subsection{Definition of the model and statement of the main result}
\label{sec:definition}
We are given a $d$-regular infinite rooted tree, $d\ge3$. The root is
denoted by $\rho$ and the level/depth of a vertex $v$ in the tree by
$|v|$. The notation $u\sim v$ denotes that $u$ and $v$ are connected by
an edge. The parent of a vertex $v$ is denoted by $\parent v$ (with the
convention $\parent\root= \root$). We write $u\le v$ if $u$ is an
ancestor of $v$ and $u<v$ if $u\le v$ but $u\ne v$. We furthermore
use the following handy notation: if $u\le v$, then $[u,v]$ denotes the
set of vertices on the path from $u$ to $v$, including $u$ and $v$. The
notation $(u,v]$, $[u,v)$ and $(u,v)$ then has obvious meaning. To each
edge $e=(\parent v,v)$, we then attach a random variable $X(e)$, such
that the collection $(X(e)$) is i.i.d. according to the law of a
random variable $X$. Here, orientation of the edges matters, and we
will set $X(v,\parent v) = -X(\parent v,v)$ for all~$v$.

In what follows, we will introduce several assumptions, \textit{which
we assume to hold throughout the paper}. We begin with the following
assumptions on the law of~$X$.

\begin{longlist}[(XM)]
\item[(XS)] The law of $X$ is of compact support, that is,
$\operatorname{esssup} |X| < \infty$.\vadjust{\goodbreak}
\item[(XR)] There exists $\beta_0>0$, such that $E[e^{\beta_0 X}
f(X)]=E[f(-X)]$ for all bounded measurable functions $f$.
\item[(XM)] $\inf_{\beta\ge0}\Lambda(\beta) > 0$, where $\Lambda
(\beta
) = \log E[e^{\beta X}] + \log(d-1)$.
\end{longlist}
Note that (XR) is equivalent to the Laplace transform $\beta\mapsto\E
[e^{\beta X}]$ being symmetric around $\beta_0/2$. In particular, the
constant $\beta_0$ is necessarily unique unless $X=0$ almost surely.

An example for a law satisfying (XS), (XR) and (XM) is the distribution
of $2Y - n$, where $Y$ follows a binomial distribution of parameters
$n$ and $p$, with $p\in(p_0,1/2)$, where $p_0 = (1-\sqrt{1-(d-1)^{-2/n}})/2$.
In this case, $\beta_0=\log\frac{1-p}p$. In general, in order to
construct a law satisfying (XS) and (XR), one can start from a
symmetric random variable $X$ taking values in a compact interval
$[-K,K]$ and define a law with Radon--Nikodym derivative proportional
to $e^{(-\beta_0/2) X}$ with respect to the law of $X$. This law will
then satisfy (XM) for $\beta_0$ small enough.

We remark that assumption (XS) seems not to be crucial, and the
argument extends
to certain distributions with non compact support, at the cost of more
complicated technical arguments. To avoid this complication we chose to present
the result under this simplifying assumption. On the other hand, assumption
(XR) is essential for our treatment, as it ensures, at $\beta_0$,
the reversibility of the
Markov chain consisting of the environment viewed from the point of
view of the particle; see Proposition~\ref{prop:rev_erg}. The reversibility
will be crucial both in the application of the Kipnis--Varadhan theory,
as well
as in the proof of
validity of the Einstein relation (one may expect a correction term
for non reversible chains).

We now define the branching random walk by
%
%
\begin{equation}
\label{eq:S} S(v) = \sum_{u\in(\root,v]} X(\parent u,u),\qquad S(
\root) = 0.
\end{equation}
Note that $X(u,v) = S(v) - S(u)$ for every two vertices $u$ and $v$
with $u\sim v$, by the above convention that $X(u,v) = -X(v,u)$. Since
$\Lambda(\beta)$ is the log-Laplace transform of this branching random
walk, it is known \cite{Biggins1977} that $\lim_{n\to\infty} \max
_{|v|=n} S(v)/n$ exists and is positive under assumption (XM). Note
further that assumptions (XM) and (XR) together imply that $\Lambda
(\beta) > 0$ for all $\beta\in\R$, such that $\lim_{n\to\infty}
\min_{|v|=n} S(v)/n$ exists and is negative \cite{Biggins1977}.

In order to define the Metropolis algorithm, we are given a function $
h\dvtx\R_+\to\R$ satisfying the following conditions [examples are
$h(x) =
\min(1,x)$ and $h(x) = x/(1+x)$]:
\begin{longlist}[(H2)]
\item[(H1)] $h$ takes values in $[0,1]$, is nondecreasing and
satisfies $h(0) = 0$ and $\lim_{x\to\infty} h(x) = 1$.
\item[(H2)] It is Lipschitz-continuous and continuously differentiable
on $(0,1)\cup(1,\infty)$.
\item[(H3)] It satisfies the functional equation $h(x) = xh(1/x)$ for
all $x\ge0$.\vadjust{\goodbreak}
\end{longlist}

For a given realization of the branching random walk and a parameter
$\beta\in\R$, the Metropolis algorithm is then the Markov chain
$(V_n)_{n\ge0}$ on the vertices of the tree with the transition
probabilities $P_{\beta}(v,w)$ given by
\begin{eqnarray*}
P_{\beta}(v,w) &=&p_{\beta}\bigl( X(v,w)\bigr)\qquad\mbox{for }
w\sim v,
\mbox{ where } p_\beta(x) = \frac{1} d h\bigl(e^{(\beta_{0}+\beta
) x}
\bigr),
\\
P_{\beta}(v,v) &= &1 - \sum_{w\sim v}
P_{\beta}(v,w).
\end{eqnarray*}
%
We denote the (annealed, i.e., averaged over the environment)
law of the Metropolis algorithm on the branching random walk tree by
$\PP_\beta$ and expectation with respect to this law by $\EE_\beta$.
Our main theorem is the following:

%
\begin{theorem}
\label{th:einstein} Set $S_n = S(V_n)$.
\begin{longlist}[(1)]
\item[(1)] The limit $\sigma^2 = \lim_{n\to\infty} S_n^2/n$
exists $\PP_0$-almost surely and is a strictly
positive and finite constant.
\item[(2)] For each $\beta\in\R$, the deterministic limit $
v_\beta=
\lim_{n\to\infty} S_n/n$ exists $\PP_\beta$-almost surely and satisfies
%
%
\begin{equation}
\label{eq-einsteinfinal} \lim_{\beta\to0} \frac{v_\beta}{\beta} =
\frac{\sigma^2} 2.
\end{equation}
\end{longlist}
\end{theorem}
We note that the existence of $v_\beta$ and the fact that it vanishes at
$\beta=0$ were already shown in \cite{Aldous1998}
(in a slightly more restrictive setup).
The main novelty in Theorem~\ref{th:einstein} is the proof of the Einstein relation
\eqref{eq-einsteinfinal}, as well as the fact that the right side is strictly
positive.
\subsection{Related works}
Our main inspiration, as noted above, is Aldous's work~\cite
{Aldous1998}. In that
paper, Aldous makes the crucial observation that a reversible invariant measure
for the environment viewed from the point of view of the particle
exists at
$\beta=0$, and derives from this that $v_0=0$, and the existence of
the limit
$\sigma^2$ under $\PP_0$; he also completely analyzes a greedy
algorithm and
formulates a series of conjectures, some answered here.
In the same paper, Aldous also refers speculatively
to \cite{Lebowitz1994} as relevant to the analysis near $\beta=0$;
indeed, the
approach of the latter to proofs of the Einstein relation forms the basis
of the current paper, as well to recent advances in the analysis of the Einstein
relation for disordered systems, as we now discuss.

The Einstein relation (ER) links the asymptotic variance of additive functionals
of (reversible) Markov chains in
equilibrium to the chains' response to small perturbations.
In a weak limit (where the time-scale is related to the strength of the
perturbation), Lebowitz and Rost \cite{Lebowitz1994} provide a general recipe
(based on the Kipnis--Varadhan theory, see
\cite{Komorowski2012}
for a comprehensive account) for the validity of a weak form
of the ER in disordered systems. For
the tagged particle in the symmetric exclusion process, the
ER was proved by Loulakis in $d\ge3$ \cite{LOU} by perturbative
methods (using transience in an essential way); this
approach
was adapted to bond diffusion in $\mathbb
Z^d$ in special environment distributions \cite{Komorowski2005a}. For mixing
dynamical random
environments with spectral gap, a full perturbation expansion was
proved in
\cite{OK1}.

Significant recent progress was achieved by \cite{Gantert2012}, where the
Lebowitz--Rost approach was combined
with good \textit{uniform in the environment}
estimates on certain regeneration times in the transient regime, that
are used
to pass from a weak ER to a full ER.
These uniform estimates are typically
not available for random walks on (random) trees, and a completely different
approach, based on explicit recursions,
was taken in \cite{BenArous2011},
where (biased)
random
walks on Galton--Watson trees were analyzed.
While we still consider walks on trees, the approach we take is closer
to that of \cite{Gantert2012}, while
replacing their uniform regeneration estimates with probabilistic estimates,
in the spirit of
\cite{Peres2008}. See also \cite{Guo2012} for another approach to
the proof of the ER in the context of balanced random walks.

\subsection{Overview of the proof and outline of the paper}
As mentioned above, the starting point is Aldous's observation
that under $\PP_0$,
the environment viewed from the point of view of the particle forms a
\textit{reversible} Markov chain. We begin by proving this (Proposition~\ref{prop:rev_erg}), and then apply the Kipnis--Varadhan theory
to deduce an invariance principle for anti-symmetric
additive functionals (Lemma~\ref{lem:antisym}). This allows us to prove the weak ER,
Theorem~\ref{th:weak_er}, following the Lebowitz--Rost recipe.

To handle the
perturbation, estimates on regeneration times and distances are crucial.
We work with \textit{level} regeneration times that are introduced
in Section~\ref{sec-reg};
these involve the random walk location $\{V_n\}$,
not
the vertices values $\{S_n\}$; of course, the latter influence the transition
probabilities of the random walk.
In order to transform the weak ER to a full ER,
we need uniform bounds on the moments of the regeneration times. These
are obtained in
Proposition~\ref{prop:tau_moments}, where it is proved that the regeneration
times exhibit uniform annealed stretched-exponential bounds. The proof has
two main steps: first, exponential moments are
proved for regeneration \textit{distances},
using in a crucial way a structure lemma of Grimmett and Kesten; see
Lemma~\ref{lem:excursion_tail}.
Then, the estimates for regeneration times are obtained, using
that the walk
must visit many well-separated fresh vertices, and between
two such visits, the walk has a large enough probability to hit distant
levels.
In proving the last statement, an argument of Aid\'{e}kon \cite{Aidekon2008}
is used; see Lemma~\ref{lem:fresh}.

\section{The weak Einstein relation}
\label{sec:weak_er}

In this section, we show that the Einstein relation holds for times of
the order of $\beta^{-2}.$ Specifically, we will prove the following result:

%
\begin{theorem}
\label{th:weak_er}
\textup{(1)} $\EE_0[S_n]=0$ for all $n\ge0$. Furthermore, the limit
$\sigma^2 = \lim_{n\to\infty} \EE_0[S_n^2]/n$ exists and is a finite,
nonnegative constant.

\textup{(2)} Set $S^\beta= \beta S_{\lfloor\beta^{-2} \rfloor}$.
Then,
\[
\EE_\beta\bigl[S^\beta\bigr] \to\frac{\sigma^2} 2 \qquad\mbox{as }\beta
\to0.
\]
\end{theorem}

The proof of Theorem~\ref{th:weak_er} uses a fairly generic and now
classical change of measure argument in the spirit of Lebowitz and Rost
\cite{Lebowitz1994}. It uses the crucial concept of the \emph
{environment seen from the particle}, which we define as follows.

Let $\Omega$ be the space of rooted $d$-regular unlabeled trees $
\omega
$ with marked edges, that is, to every two vertices $u$ and $v$ with
$u\sim v$ we associate a real number $X_\omega(u,v)$ with $X_\omega
(v,u) = -X_\omega(u,v)$. Note that ``unlabeled'' means that we do not
distinguish between the neighbors of a vertex,\setcounter{footnote}{2}\footnote{There are
several ways how to render this formal, one of which consists of first
defining the space $\widetilde\Omega$ of labeled rooted $d$-regular
trees with marks, which is homeomorphic to $\R^\N$. The space $\Omega$
is then defined as the quotient space with respect to the group of
graph automorphisms fixing the root. It is endowed with the Borel
$\sigma$-algebra induced by the quotient topology.} this will be
crucial for what follows. The root of every tree $\omega$ is denoted
by $\root$. For every vertex $v$, we then define the shift operator
$\theta_v\dvtx\Omega\to\Omega$, which yields the tree $\omega$
``seen from
the vertex $v$.'' A bit more formally,\footnote{In order to render this
completely formal, one can first define the shift operator on the
auxiliary space $\widetilde\Omega$ (see above) and then show that it
induces a well-defined operator on $\Omega$.} if for a vertex $u$ in
$\omega$ we denote by $\theta^{-1}_vu$ its corresponding vertex in
$\theta_v\omega$, then
$X_{\theta_v\omega}(\theta^{-1}_vu,\theta^{-1}_vw) = X_\omega(u,w)$.
In particular, if $v\sim\rho$, then the mark of the edge $(\rho,v)$
``changes its sign upon passing from $\omega$ to $\theta_v \omega$.''

Define $S_\omega(v)$ for each vertex $v$ analogously to \eqref{eq:S}.
For $\omega\in\Omega$, we then define the operators $\Lc_\omega$ and
$\Lt_\omega$ acting on functions $f\dvtx\R\to\R$ by
\[
\Lc_\omega f = \sum_{v\sim\root} f
\bigl(S_{\omega}(v)\bigr) = \sum_{v\sim
\root} f
\bigl(X_{\omega}(\root,v)\bigr)\quad \mbox{and}\quad\Lt_\omega f =
\Lc_\omega(pf),
\]
where $p(x) = p_0(x) = d^{-1} h(e^{\beta_0 x})$; see Section~\ref
{sec:definition}.

Let $\P$ be the law on $\Omega$ under which all edges $X_\omega
(\parent
u,u)$ are i.i.d. according to the law of $X$; see Section~\ref
{sec:definition}. We denote by $\E$ the expectation with respect to
$\P$.
Furthermore, let $\PP_{\beta}$ be the law of the Metropolis algorithm
$(V_n)_{n\ge0}$ with transition probabilities $P_\beta$ defined in
Section~\ref{sec:definition} and the underlying tree $\omega_0$
distributed according to $\P$. Expectation w.r.t. $\PP_\beta$ is
denoted by $\EE_\beta$, and we also set $\PP= \PP_0$ and $\EE= \EE
_0$. Setting $\omega_n = \theta_{V_n}\omega_0$ then defines a Markov
chain $(\omega_n)_{n\ge0}$ on the space $\Omega$, which jumps from
$\omega$ to $\theta_v\omega$ with probability $p_{\beta}(X(\root,v))$
for every $v\sim\root$. Let $(\F_n)_{n\ge0}$ be the natural filtration
of $(\omega_n)_{n\ge0}$, augmented by sets of zero measure. The process
$S_n = S_{\omega_0}(V_n)$ is then adapted to $(\F_n)_{n\ge0}$, since it
can be almost surely reconstructed from $\omega_0,\ldots,\omega_n$.

The following result was already observed by Aldous \cite{Aldous1998},
who had a more complicated proof for it.
%
%
\begin{proposition}
\label{prop:rev_erg}
The process $(\omega_n)_{n\ge0}$ is reversible and ergodic under~$\PP$.
\end{proposition}
\begin{pf}
In order to show reversibility, since $(\omega_n)_{n\ge0}$ is a Markov
process, we only have to show that $\EE[F(\omega_0,\omega_1)] = \EE
[F(\omega_1,\omega_0)]$ for every bounded (Borel) measurable functional
$F\dvtx\Omega^2\to\R$. For this, it is obviously enough to show that
$\EE
[(F(\omega_0,\omega_1) - F(\omega_1,\omega_0))\Ind_{\omega_0 \ne
\omega
_1}] = 0$. Now we have
\begin{eqnarray*}
\EE\bigl[F(\omega_0,\omega_1)\Ind_{\omega_0 \ne\omega_1}
\bigr] &=& \sum_{v\sim
\rho}\EE\bigl[p\bigl(X_{\omega_0}(
\rho,v)\bigr)F(\omega_0,\theta_v\omega_0)
\bigr]
\\
&=& \sum_{v\sim\rho}\E\bigl[p\bigl(-X_{\omega}(
\rho,v)\bigr)e^{\beta_0X_{\omega
}(\rho
,v)}F(\omega,\theta_v\omega)\bigr],
\end{eqnarray*}
where the last equality follows from assumption (H3). Conditioned
on\break
$X_{\omega}(\rho,v)$, the environment $\theta_v\omega$ is distributed
as $\omega$ but with one edge pointing away from the root bearing the
value $-X_\omega(\rho,v)$ (remember that the vertices are unlabeled,
such that ``it can be any one of them,'' which amounts to saying that
``we do not know where we came from''). By assumption (XR), the
right-hand side of the last equation is therefore equal to
\[
\sum_{v\sim\rho}\E\bigl[p\bigl(X_{\omega}(\rho,v)
\bigr)F(\theta_v\omega,\omega)\bigr] = \EE\bigl[F(
\omega_1,\omega_0)\Ind_{\omega_0 \ne\omega_1}\bigr],
\]
which finishes the proof of the reversibility.

Ergodicity follows from a classical ellipticity
argument which we recall (see also \cite{Zeitouni2004},
Corollary~2.1.25, for a similar argument): Let $Q$ be a
stationary probability measure of the Markov chain $(\omega
_{n})_{n\ge
0}$ with $Q\ll\P$. We wish to show that $\P\ll Q$, which will imply
ergodicity since ergodic measures are the extremal points in the convex
set of stationary probabilities. Define the event $E=\{\dd Q/\dd\P=
0\}$. By invariance, $E_Q[P\Ind_E]= E_{Q}[\Ind_E]=0$, where $P$ is the
transition kernel of the Markov chain. This further implies $\Ind_E
\ge
P\Ind_E$, $\P$-almost surely. Since the transition probabilities are
strictly positive, we then have $\Ind_E(\omega) \ge\max_{v\sim\rho}
\Ind_E(\theta_v\omega)$, because $\Ind_E(\omega)$ takes values in
$\{
0,1\}$. Fixing an infinite ray $\rho= v_0,v_1,v_2,\ldots,$ we then get
by iteration of the previous inequality that $\Ind_E(\omega)\ge\Ind
_E(\theta_{v_i}\omega)$ for every~$i$, whence $\Ind_E \ge n^{-1}\sum
_{i=1}^n \Ind_E(\theta_{v_i}\omega)$ for every $n$. But since $\P
$ is
a pro\-duct measure and therefore ergodic with respect to the shift
along the ray, Birkhoff's ergodic theorem now gives $\Ind_E\ge\P(E)$,
$\P$-almost surely, which implies $\P(E)\in\{0,1\}.$ But $Q\ll\P
$ by
hypothesis, whence $\P(E)=0$. This finishes the proof.
\end{pf}

We recall the following basic fact about reversible processes.
%
%
\begin{lemma}
\label{lem:rev_functional}
For any bounded measurable functionals $F$ and $G$ and every $n\ge0$,
we have
\[
\EE\bigl[F(\omega_0,\ldots,\omega_n)G(
\omega_n)\bigr] = \EE\bigl[F(\omega_n,\ldots,
\omega_0)G(\omega_0)\bigr].
\]
\end{lemma}

We will need the following result about anti-symmetric additive
functionals of reversible ergodic Markov processes. It is implicit in
the proof of Theorem~2.1 in~\cite{DeMasi1989} and relies on a
celebrated result from \cite{Kipnis1986a}; see also Chapters~1 and 2 in
\cite{Komorowski2012} for a comprehensive account of the theory.
%
%
\begin{lemma}
\label{lem:antisym}
Let $F\dvtx\Omega^2\to\R$ be an anti-symmetric measurable
functional, that
is, $F(\omega,\omega') = -F(\omega',\omega)$ for all $\omega
,\omega'\in
\Omega$, with $\EE[F(\omega_0,\omega_1)^2] < \infty$. Define a sequence
of random variables by $\mathfrak S_n = \sum_{k=1}^n F(\omega
_{k-1},\omega_k)$ for all $n>0$. Then:
\begin{longlist}[(1)]
\item[(1)] The time-variance $\sigma^2 = \lim_{n\to\infty} \frac
{1} n
\EE
[\mathfrak S_n^2]$ exists and is finite.
\item[(2)] There exists a square integrable martingale $M_n$ with
stationary ergodic increments, such that $\frac{1} {\sqrt n}(\mathfrak
S_n - M_n)$ converges to 0 in $L^2$. In particular, $\frac{1} {\sqrt n}
\mathfrak S_n$ converges in law to a centered Gaussian variable with
variance~$\sigma^2$.
\end{longlist}
\end{lemma}

The following lemma makes precise an expansion of $p_\beta$ around $p$
for small~$\beta$. It easily follows from assumptions (H1)--(H3).
%
%
\begin{lemma}
\label{lem:p_expansion}
There exist measurable functions $q_\beta(x)$ and $q(x)$, such that:

\begin{longlist}[(1)]
\item[(1)] $p_\beta(x)=p(x)\exp(\beta q_\beta(x))$;
\item[(2)] $q_\beta(x)$ is uniformly bounded for all $x\in\R$ and small
enough $\beta$, and $q_\beta(x)\to q(x)$ as $\beta\to0$, for all
$x\in
\R$;
\item[(3)] $q(0) = q_\beta(0) = 0$;
\item[(4)] there exists a constant $c>0$, such that for small enough
$\beta$, we have $|\Lt_{\omega}(e^{\beta q_\beta}-1)| \le\beta
c(1-\Lc
_{\omega} p)$ and $|\Lt_\omega q| \le c(1-\Lc_{\omega} p)$.
\end{longlist}
\end{lemma}

We are now ready for the proof of Theorem~\ref{th:weak_er}.
\begin{pf*}{Proof of Theorem~\ref{th:weak_er}}
\textit{First part}: We first note that there exists a measurable
functional $F\dvtx\Omega^2\to\R$, such that $S_1 = F(\omega
_0,\omega_1)$,
$\PP$-almost surely. This follows from the fact that $\P$-almost
surely, the $d$ shifted environments $\theta_v \omega$, $v\sim\root$
are all different, otherwise there would be at least two identical
subtrees of the vertices in the second generation which is an event of
probability zero (except if $X=0$ almost surely, in which case the
lemma is trivial). By the definition of $S_1$, we can furthermore
choose $F$ to be anti-symmetric in the sense of Lemma~\ref
{lem:antisym}. In particular, $\EE[S_1] = -\EE[S_1] = 0$ by
Lemma~\ref
{lem:rev_functional}, whence $\EE[S_n]=0$ for all $n\ge0$. The second
statement follows from the first part of Lemma~\ref{lem:antisym}.

\textit{Second part}: We will use a change of measure argument as in
\cite{Lebowitz1994}. The basic idea is to write the Radon--Nikodym
derivative of $\PP_\beta$ with respect to $\PP$ as an exponential
martingale of the form $\exp(Z^\beta_n - \frac{1} 2 A^\beta_n)$ for some
martingale $(Z^\beta_n)_{n\ge0}$ and to show that the pair $(\beta
S_{\lfloor\beta^{-2}\rfloor},Z^\beta_{\lfloor\beta^{-2}\rfloor})$
converges in law under $\PP$ to a centered Gaussian vector $(G_S,G_Z)$
with covariance $\EE[G_SG_Z] = \frac{1} 2 \EE[G_S^2] = \frac{1} 2
\sigma
^2$. The theorem then follows from a standard change of measure
argument for Gaussian variables. Here are the details:

\emph{Step} 0: Set $\Delta S_n = S_{n+1}-S_n$ for all $n\ge0$. The
Radon--Nikodym derivative of $\PP_\beta$ with respect to $\PP$ is
given by
%
%
\begin{equation}\qquad
\label{eq:2} \log\frac{\dd\PP_\beta}{\dd\PP}\bigg |_{\F_n} = \sum
_{k=0}^{n-1} \biggl[\log\frac{p_\beta(\Delta S_k)}{p(\Delta S_k)} +
\log
\biggl(\frac
{1-\Lc
_{\omega_k}p_\beta}{1-\Lc_{\omega_k} p} \biggr)\Ind_{\omega_k =
\omega
_{k+1}} \biggr]=: Y_n.
\end{equation}
Note that if $1-\Lc_{\omega_k} p=0$, then $\omega_k \ne\omega_{k+1}$
with probability 1, such that the second summand is well defined. Note
also that here and in what follows, the empty sum always has value 0.
If $q_\beta$ and $q$ are the functions from Lemma~\ref
{lem:p_expansion}, we can write the process $Y_n$ as
%
%
\begin{equation}
\label{eq:Y} Y_n = \beta\sum_{k=0}^{n-1}
q_\beta(\Delta S_k) + \sum_{k=0}^{n-1}
\log\biggl(1-\frac{\Lt_{\omega_k}(e^{\beta q_\beta}-1)}{1-\Lc
_{\omega_k}
p} \biggr)\Ind_{\omega_k = \omega_{k+1}}.
\end{equation}
We then define the process $(A_n)_{n\ge0}$ by (we suppress the
dependence on $\beta$ from the notation)
\begin{eqnarray*}
A_n &=& -2 \sum_{k=1}^n
\EE[Y_k-Y_{k-1}|\F_{k-1}] \\
&= &-2\sum
_{k=0}^{n-1} \biggl[\beta\Lt_{\omega_k}q_\beta+
(1-\Lc_{\omega_k}p)\log\biggl(1-\frac
{\Lt_{\omega_k}(e^{\beta q_\beta}-1)}{1-\Lc_{\omega_k}p} \biggr)
\biggr].
\end{eqnarray*}
(Consistent with the definitions we have that $Y_0=A_0=0$.)
Note that by Lemma~\ref{lem:p_expansion}, we have for small $\beta$,
%
%
\begin{eqnarray}
\label{eq:An} A_{n+1} - A_n &=& 2 \bigl(
\Lt_{\omega_n}\bigl(e^{\beta q_\beta}-1\bigr) \beta\Lt_{\omega
_n}q_\beta
\bigr) - \frac{ (\Lt_{\omega_n}
(e^{\beta q_\beta}-1) )^2}{1-\Lc_{\omega_n}p} + O\bigl(\beta
^3\bigr)
\nonumber
\\[-8pt]
\\[-8pt]
\nonumber
&=& \beta^2 \bigl(\Lt_{\omega_n} q^2+(
\Lt_{\omega_n} q)^2/(1-\Lc_{\omega
_n} p) + o(1) \bigr).
\end{eqnarray}
We further define the process $(Z_n)_{n\ge0}$ (again suppressing the
dependence on $\beta$) by
\[
Z_n = Y_n - \sum_{k=1}^n
\EE[Y_k-Y_{k-1}|\F_{k-1}] = Y_n +
\frac{1} 2 A_n,
\]
such that $(Z_n)_{n\ge0}$ is a martingale under $\PP$ with respect to
the filtration $(\F_n)_{n\ge0}$.

\emph{Step} 1: We wish to show that the random variable $Z^\beta=
Z_{\lfloor\beta^{-2} \rfloor}$ converges in law under $\PP$ to a
centered Gaussian variable with variance $\sigma_Z^2 < \infty$.
Define the $\PP$-martingale $(M_n)_{n\ge0}$ by
\[
M_n =\sum_{k=0}^{n-1} \biggl(q(
\Delta S_k) - \frac{\Lt_{\omega_k}
q}{1-\Lc_{\omega_k} p}\Ind_{\omega_k = \omega_{k+1}} \biggr).
\]
By \eqref{eq:An} and Lemma~\ref{lem:p_expansion}, we then have
%
%
\begin{eqnarray}
\label{eq:Zn_Mn} \EE\bigl[(Z_n - \beta M_n)^2
\bigr] &=& \EE\Biggl[\sum_{k=0}^{n-1}
\bigl[Y_{n+1} - Y_n - \beta(M_{n+1} -
M_n) + O\bigl(\beta^2\bigr)\bigr]^2 \Biggr]
\nonumber
\\[-8pt]
\\[-8pt]
\nonumber
& =&
o\bigl(\beta^2 n\bigr),
\end{eqnarray}
whence $\EE[(Z^\beta- \beta M_{\lfloor\beta^{-2} \rfloor})^2] \to0$
as $\beta\to0$.

Note that by the fourth point of Lemma~\ref{lem:p_expansion}, $M_n$ is
square-integrable and by Proposition~\ref{prop:rev_erg}, the sequence
of its increments $(M_{n+1}-M_n)_{n\ge0}$ is stationary and ergodic.\vadjust{\goodbreak}
By the martingale central limit theorem for stationary ergodic
sequences (see, e.g., \cite{Durrett1996}, Theorem~7.7.5),\vspace*{1pt} the sequence
$(M_n/\sqrt n)_{n\ge0}$ then converges in law under $\PP$ to a centered
Gaussian variable with variance $\sigma_Z^2 = \EE[M_1^2] = \EE[\Lt
_\omega q^2+(\Lt_\omega q)^2/(1-\Lc_\omega p)]$. Together with \eqref
{eq:Zn_Mn}, this proves the above-mentioned convergence of $Z^\beta$.

\emph{Step} 2: We wish to show that the random variable $A^\beta=
A_{\lfloor\beta^{-2} \rfloor}$ converges in probability to $\sigma
_Z^2$ under the law $\PP$. Define the process $A'_n$ by
$
A'_n = \sum_{k=0}^{n-1} [\Lt_{\omega_k} q^2+(\Lt_{\omega_k}
q)^2/(1-\Lc
_{\omega_k} p)].
$ By Proposition~\ref{prop:rev_erg} and the ergodic theorem, the
sequence $(A'_n/n)_{n\ge0}$ converges $\PP$-almost surely to $\EE[A'_1]
= \sigma_Z^2$. Together with \eqref{eq:An}, this yields the
above-mentioned convergence of $A^\beta$ as $\beta\to0$.

\emph{Step} 3: Recall the definition $S^\beta= \beta S_{\lfloor\beta
^{-2}\rfloor}$. We wish to show that the pair $(S^\beta,Z^\beta)$
converges in law under $\PP$ to a centered Gaussian vector $(G_S,G_Z)$
with covariance $\EE[G_SG_Z] = \frac{1} 2 \EE[G_S^2] = \frac{1} 2
\sigma
^2$ (with $\sigma^2$ from Lemma~\ref{lem:antisym}). By \eqref
{eq:Zn_Mn}, it is enough to show that this convergence holds for the
pair $(S_n/\sqrt n,\break M_n/\sqrt n)$ as $n\to\infty$. By Lemma~\ref
{lem:antisym}, every linear combination $a S_n +b M_n$ is the sum of a
square-integrable martingale with stationary and ergodic increments and
a process $R_n$ with $\frac{1} {\sqrt n} R_n \to0$ in $L^2$, as $n\to
\infty$. Again by the martingale CLT for stationary, ergodic sequences
(\cite{Durrett1996}, Theorem~7.7.5), the pair $(S_n/\sqrt n,M_n/\sqrt
n)$ then converges in law as $n\to\infty$ to a centered Gaussian vector
$(G_S,G_Z)$ with $\EE[G_S^2] =\sigma^2$, $\EE[G_Z^2] = \sigma_Z^2$ and
$\EE[G_SG_Z] = \lim_{n\to\infty} \EE[S_nM_n]/n = \lim_{\beta\to
0}\EE
[S^\beta Z^\beta]$.

It remains to show that $\lim_{\beta\to0}\EE[S^\beta Z^\beta]=
\sigma
^2/2$. In order to prove this, recall the definition of $\Delta
S_{n}=S_{n+1}-S_n$ and define $\Delta Z_{n}=Z_{n+1}-Z_n$. We have
for every $n\ge0$,
%
%
\begin{equation}
\label{eq:3} \EE\bigl[S_{n+1}^2-S_n^2
\bigr] = \EE\bigl[(\Delta S_n)^2 + 2S_n\Delta
S_n\bigr] = \EE\bigl[\Lt_{\omega_n}\bigl(x^2\bigr)
+ 2S_n \Lt_{\omega_n} x\bigr],
\end{equation}
where $x$ is the identity function. By Lemma~\ref{lem:rev_functional},
we have $\EE[S_n\Lt_{\omega_n} x] =  \EE[(-S_n)\Lt_{\omega_0} x]$,
whence, summing \eqref{eq:3} over $n,$ we get,
\[
\EE\bigl[S_n^2\bigr] = \sum
_{k=1}^n \EE\Biggl[\Lt_{\omega_0}
\bigl(x^2\bigr) - 2 \Lt_{\omega
_0} x \times\sum
_{j=0}^{k-2} \Delta S_j \Biggr].
\]
Furthermore, since $(Z_n)_{n\ge0}$ is a martingale, we have
%
%
\begin{eqnarray}
\label{eq:5} \EE[S_{n+1}Z_{n+1} - S_nZ_n]
&=& \EE[Z_{n+1}\Delta S_n]
\nonumber
\\[-8pt]
\\[-8pt]
\nonumber
& =& \EE\bigl[\beta\Delta
S_n \bigl(q_{\beta}(\Delta S_n)-
\Lt_{\omega_n} q_\beta\bigr)+(\Lt_{\omega_n}
x)Z_n \bigr],
\end{eqnarray}
where we made use of the fact that $\Delta S_k = 0$ on the event that
$\omega_k = \omega_{k+1}$. Applying Proposition~\ref{prop:rev_erg} to
the term $\EE[(\Lt_{\omega_n} x)Z_n]$, we see that the terms
corresponding to the second summand in the brackets of \eqref{eq:2}
cancel. Summing \eqref{eq:5} over $n$, this yields
\begin{eqnarray*}
\EE[S_nZ_n] &=& \beta\sum_{k=1}^n
\EE\Biggl[\Lt_{\omega_0}(xq_{\beta
}) - \Lt_{\omega_0} x \times
\sum_{j=0}^{k-2} \bigl(q_\beta(
\Delta S_j) - q_\beta(-\Delta S_j) \bigr)
\Biggr] \\
&&{}+ \EE[S_{n-1}\Lt_{\omega_0} q_\beta].
\end{eqnarray*}
Now, by assumption (H3) we have $q_\beta(x) - q_\beta(-x) = x$ for
every $x$ (this is the critical point!). Moreover, by reversibility, we
have $\E[\Lt_\omega f] = \E[\Lt_\omega\widebar f]$ for every
function $f$, where $\widebar f(x) = f(-x)$. This yields
\[
\E\bigl[\Lt_{\omega}(xq_{\beta}) \bigr] = \tfrac1 2 \times\E
\bigl[\Lt_{\omega}(xq_{\beta}-x\widebar{q}_\beta) \bigr] =
\tfrac1 2 \times\E\bigl[\Lt_{\omega}\bigl(x^2\bigr) \bigr].
\]
%
Altogether, the previous equations now yield
\[
\EE\bigl[S^\beta Z^\beta\bigr] = \EE\bigl[
\bigl(S^\beta\bigr)^2\bigr]/2+ \beta\EE[S_{\lfloor\beta
^{-2}-1\rfloor}
\Lt_{\omega_0}q_\beta].
\]
Convergence of the first summand has been established above, and the
second tends to 0 by Lemma~\ref{lem:antisym} and the Cauchy--Schwarz
inequality. Hence, we obtain $\lim_{\beta\to0}\EE[S^\beta Z^\beta]=
\sigma^2/2$ as claimed.

\emph{Step} 4: We claim that $\EE_\beta[S^\beta] \to\EE[G_S\exp
(G_Z -
\frac{1} 2 \sigma_Z^2)]$, as $\beta\to0$. Since $(G_S,G_Z)$ is a
centered Gaussian vector with $\EE[G_SG_Z] = \frac{1} 2 \EE[G_S^2] =
\sigma^2/2$ and $\EE[G_Z^2]=\sigma_Z^2$, this will finish the proof of
the theorem. By \eqref{eq:2}, $\EE_\beta[S^\beta] = \EE[S^\beta
\exp
(Z^\beta- \frac{1} 2 A^\beta)]$ and by the convergences in law
established above, it suffices to show that this last expression is
uniformly integrable. Now, since $Z_n$ and $\exp(Z_n -\frac{1} 2 A_n)$
are martingales, $A_n$ is a submartingale. $A_n$ being $\F
_{n-1}$-measurable, it is therefore increasing in $n$. It then remains
to show that $S^\beta\exp(Z^\beta)$ is uniformly integrable. By the
fourth point of Lemma~\ref{lem:p_expansion}, $Z_n$ is a martingale with
bounded increments for $\beta$ small enough. Azuma's inequality \cite
{Azuma1967} then implies that all exponential moments of $Z^\beta$ are
uniformly bounded in $\beta$, for small enough $\beta$. Furthermore,
$\EE[(S^\beta)^2]$ is uniformly bounded by the first part of this
theorem. H\"older's inequality then yields
uniform boundedness of $\EE[(S^\beta\exp(Z^\beta))^c]$ for some
constant $c>1$, which finishes the proof.
\end{pf*}

\section{Estimates on the branching random walk}
\label{sec:brw}
In this section, we establish an estimate for the branching random walk
(Lemma~\ref{lem:brw_all_good} below). We recall that for two vertices
$u,v$ with $u\le v$, we denote by $[u,v]$ the set vertices on the path
connecting $u$ and $v$. Similarly, if $n,m\in\N$, then we define
$[n,m]$ to be the set of vertices between levels/depths $n$ and $m$.
More generally, for a vertex $u$, we let $[n,m]_u$ denote the set of
vertices between levels/depths $n$ and $m$ in the subtree rooted at $u$
(which means that these vertices are between levels $n+|u|$ and $m+|u|$
in the original tree), such that $[m,n] = [m,n]_\root$. Finally, we
write $[n]$ for $[n,n]$ and $[n]_u$ for $[n,n]_u$.

%
\begin{lemma}
\label{lem:brw_root_great}
There exist $c\in(0,\infty)$ and $b > 1$, such that for large $L$,
\[
\P\Bigl(\forall v\in[L]\dvtx\max_{w\in[\root,v]} \bigl|S(w)\bigr| > c\log L
\Bigr) \le e^{-L^b}.
\]
\end{lemma}

We will first establish the following intermediate bound:
%
%
\begin{lemma}
\label{lem:brw_root_good}
There exist constants $C_1,C_2>0$ such that
for all large $L$,
%
\[
\P\Bigl(\exists v\in[L]\dvtx\max_{w\in[\root,v]} \bigl|S(w)\bigr| \le
C_1\Bigr) > C_2.
\]
\end{lemma}
\begin{pf}
The proof is a standard first and second moment calculation.
Fix $C_1>0$ large. Let
$(S_n)_{n\ge0}$ be a random walk starting at $0$ with
steps distributed according to the law of $X$. For $n\ge0$ and
$x\in[-C_1,C_1]$, define the event
$B_n^{(x)} = \{\forall k\le n\dvtx|S_k+x| \le C_1\}$ and set $B_n = B_n^{(0)}$.
By assumption (XM) and standard large and small deviations estimates,
there exists $c_0<d-1$, such that for $C_1$ large enough,
%
%
\begin{equation}
\label{eq-lsdev} \exists L_0\in\N\ \forall L>L_0\ \forall
n\ge0\qquad \P(B_n) \ge c_0^{-n}.
\end{equation}
Indeed, this is obtained, for example, by combining
the change of measure in the Mogulskii--Varadhan theorem
\cite{DeZe}, Theorem~5.1.2, with
the fact that a centered random walk with bounded i.i.d. increments
stays in a tube of width $a$ for time $n$ with probability
at least $e^{-Cn/a^2}$ for all $n$ and $a>a_0$ and some constant $C>0$;
for finer estimates see, for example, \cite{Mogulskii}.

In the sequel, we fix such a $C_1$ once and for all.
By an argument similar to the above,
there exists a constant $C_1'$ depending on $C_1$ only such that
%
\begin{equation}
\label{eq:Bnx} \bigl(C_1'\bigr)^{-1}\sup
_{x\in[-C_1,C_1]} \P\bigl(B^{(x)}_n\bigr) \le
\P(B_n) \le C_1' \inf_{x\in[-C_1/2,C_1/2]}
\P\bigl(B^{(x)}_n\bigr).
\end{equation}
To see \eqref{eq:Bnx}, note from the above that $P(B_{n-C})/P(B_n)$ is
bounded
by a constant depending on $C$ only, uniformly in $n>n_0(C)$,
and then couple the walk started at~$x$ with the walk started at $0$ by time
$C$, with a fixed positive probability.

The second inequality in \eqref{eq:Bnx} yields the existence of
a constant $C_1''$ (depending on $C_1$) so that for every $k\le L$,
%
%
\begin{equation}
\label{eq:BLk} \P(B_{L-k})\P(B_k) \le
C_1''\P(B_L),
\end{equation}
because conditioned on $B_k$, the probability that
$|S_k| \le C_1/2$ is bounded from below by a strictly positive constant
uniformly in $k$.

For $v\in[L]$, let
\[
A_v=\mathbf{1}_{\{\max_{w\in[\root,v]} |S(w)| \le C_1\}}.
\]
%
Further let $A = \sum_{v\in[L]} A_v$. Then,
\[
\E[A] = (d-1)^L \P\bigl(\forall n \le L\dvtx|S_n| \le
C_1\bigr) = (d-1)^L \P(B_L).
\]
As for the second moment, denote by $u\wedge v$ the most recent common
ancestor of $u$ and $v$. We then have for large $L$,
\begin{eqnarray*}
\E\bigl[A^2\bigr] &=& \sum_{u,v\in[L]}
\E[A_vA_u] \le\sum_{u,v\in[L]}
\P(B_L) \sup_{x\in[-C_1,C_1]}\P\bigl(B^{(x)}_{L-|u\wedge v|}
\bigr)
\\
& \le& C_1' C_1''
\P(B_L)^2 \sum_{u,v\in[L]}
\P(B_{|u\wedge v|})^{-1},
\end{eqnarray*}
where the last inequality follows
from \eqref{eq:Bnx} and \eqref{eq:BLk}. Equation \eqref{eq-lsdev}
now yields
\[
\sum_{u,v\in[L]}\P(B_{|u\wedge v|})^{-1} \le
\sum_{u,v\in[L]} c_0^{|u\wedge v|} \le C
(d-1)^{2L}
\]
for some $C>0$.
The lemma now follows from the previous three inequalities together
with the Paley--Zygmund bound $\P(A>0) \ge{\E[A]^{2}}/{\E[A^2]}$.
\end{pf}

\begin{pf*}{Proof of Lemma~\ref{lem:brw_root_great}}
Let $c>0$, and set $H = \lceil c \log L\rceil$. Let $C_1$ be as in
Lemma~\ref{lem:brw_root_good}.
Let $g = \operatorname{esssup} |X|$, which is finite by assumption (XS).
The branching random walks spawned by the vertices at level $H$
being independent, we have
\begin{eqnarray*}
&&\P\Bigl(\forall v\in[L]\dvtx\max_{w\in[\root,v]} \bigl|S(w)\bigr| >
gH+C_1\Bigr) \\
&&\qquad\le\P\Bigl(\forall v\in[L-H]\dvtx\max
_{w\in[\root,v]}\bigl |S(w)\bigr| > C_1\Bigr)^{(d-1)^H}.
\end{eqnarray*}
%
The lemma now follows from the last inequality together with
Lemma~\ref{lem:brw_root_good},
by choosing $c$ large enough.
\end{pf*}
%
%
\begin{lemma}
\label{lem:brw_all_good}
There exist $c\in(0,\infty)$ and $b > 1$, such that for large $L$,
\[
\P\Bigl(\exists u\in[0,L]\ \forall v\in[L]\mbox{ with } u\le v\dvtx
\max
_{ w\in[u,v]} \bigl|S(w)-S(u)\bigr| > c\log L\Bigr) \le e^{-L^b}.
\]
\end{lemma}

\begin{pf}
Let $c$ be as in the statement of Lemma~\ref{lem:brw_root_great}. We
say that a vertex $u$ is $H$-bad if for all $v\in[H]_u$ there exists
$w\in[u,v]$, such that $|S(w)-S(u)| > c\log L$. Note that if $u$ is
$H$-bad, then it is $K$-bad for every $K>H$. A simple union bound gives
\begin{eqnarray*}
\P\bigl(\exists u\in[0,L]\dvtx\mbox{ $u$ is $\bigl(L-|u|\bigr)$-bad}\bigr)
&\le&\P\bigl(
\exists u\in[0,L]\dvtx\mbox{ $u$ is $L$-bad}\bigr)\\
 &\le&(d-1)^{L+1}
\P(
\root\mbox{ is $L$-bad}).
\end{eqnarray*}
The statement then follows from Lemma~\ref{lem:brw_root_great}.
\end{pf}

\section{Regeneration times}
\label{sec-reg}

In this section, we establish a regeneration structure for the
Metropolis algorithm, which will permit us to prove Theorem~\ref
{th:einstein} from the previously established Theorem~\ref{th:weak_er}.
Recall the definition of the Metropolis algorithm $(V_n)$ from
Section~\ref{sec:definition}, which depends on a parameter $\beta\in
\R
$. Define the \emph{level regeneration} times $(\tau_n)_{n\ge0}$ by
$\tau_0 = 0$ and $\tau_{n+1}$ to be the first time after $\tau_n$ where
the chain $(V_n)_{n\ge0}$ hits a level $L$ for the first time, then
immediately jumps to level $L+1$ and never gets back to level $L$ again.

As in Sections~\ref{sec:weak_er} and~\ref{sec:brw}, we denote the law
of the branching random walk by~$\P$, which is a law on $\Omega$. We
further denote the (quenched) law of the Metropolis algorithm
$(V_n)_{n\ge0}$ started from the vertex $v$ and given the branching
random walk $\omega$ (the \emph{environment}) by $P_{\omega,\beta}^v$.
The annealed law is denoted by $\PP_\beta^v(\dd\omega,\dd V) = \P
(\dd
\omega) P_{\omega,\beta}^v(\dd V)$. We also set $\PP_\beta= \PP
_\beta
^\root$, and note that this agrees with earlier notation. Our goal is
to show:

%
\begin{proposition}
\label{prop:tau_moments}
For each $K>0$, there exists $a=a(K)>0$ and $n_a=n_a(K)>0$
such that for all $n>n_a$ and $|\beta|\leq K$,
$\PP_\beta(\tau_1 > n) \le e^{-n^a}$ and $\PP_\beta(\tau_2-\tau
_1 > n)
\le e^{-n^a}$.
\end{proposition}
The main point in Proposition~\ref{prop:tau_moments} is in
\textit{uniformity} (in $\beta$)
of the tail bounds for the regeneration times. This uniformity is in
sharp\vadjust{\goodbreak} contrast to other settings discussed
in the literature, where the regeneration times usually
blow up when the parameter approaches the critical
value \cite{BenArous2011,Gantert2011a}. We remark that we
actually only need that $\EE_\beta[\tau_1^k]$ is uniformly bounded
for $\beta$ in a neighborhood of $0$,
for some $k>2$.

In order to prove Proposition~\ref{prop:tau_moments}, we will make use
of the relation between the Markov chain $(V_n)_{n\ge0}$ and
electrical networks \cite{Lyons1992}: Let $N(v)$ be the set of
neighbors of $v$ including $v$. For $w\in N(v)$, set
\begin{eqnarray*}
Q(v,w) &=& \frac{P_\beta(v,w)}{P_\beta(v,\parent v)},\qquad Q(v) = \sum
_{w\in N(v)}
Q(v,w),
\\
C(v,w) &=& Q(v,w) \prod_{u\le v} Q(\parent u,u),\qquad C(v) =
\sum_{w\in N(v)} C(v,w).
\end{eqnarray*}
One checks that for every $w\in N(v)$, $C(v,w) = C(w,v)$ and that
\[
C(v,w)/C(v) = Q(v,w)/Q(v)= P_\beta(v,w),
\]
whence the Markov chain $(V_n)_{n\ge0}$ has an interpretation as the
random walk on the rooted $d$-regular tree with loops, induced by the
edge conductances $C(v,w)$. By assumption (H3), one has for $u\le v$,
%
%
\begin{equation}
\label{eq:conductance_ratio} \frac{C(\parent v,v)}{C(\parent u,u)} =
\frac{h(e^{(\beta_0+\beta)
X(v)})}{h(e^{(\beta_0+\beta) X(u)})} e^{(\beta_0+\beta)(S(\parent v)-S(u))}.
\end{equation}
By assumptions (XS) and (H1)--(H3), this implies the existence of a
constant $c>0$, such that
%
%
\begin{equation}
\label{eq:conductance_ratio_est} c e^{(\beta_0+\beta)(S(\parent
v)-S(u))}< \frac{C(\parent
v,v)}{C(\parent u,u)} < c^{-1}
e^{(\beta_0+\beta)(S(\parent v)-S(u))}.
\end{equation}

Define $T_L$ to be the first strictly positive time the chain
$(V_n)_{n\ge0}$ hits the level $L$. Furthermore, denote by $T_u$ and
$T_u^*$, respectively, the first nonnegative and strictly positive
times the chain hits a vertex $u$.

The first lemma gives a uniform bound on the annealed probability that
the Metropolis algorithm started from a vertex $v$ escapes to infinity
without coming back to its parent $\parent v$. It was essentially
already observed by Aldous~\cite{Aldous1998}, Lemma~8.
%
%
\begin{lemma}
\label{lem:transience}
For each $K>0$, there exists $c=c(K)
>0$, such that for each vertex $v\ne\rho$ and for all $|\beta|\le K$,
we have
\[
\E\bigl[P^v_{\omega,\beta}(T_{\parent v} = \infty)\bigr] > c.
\]
\end{lemma}
\begin{pf}
Fix $v\ne\rho$ and define $f(\beta):= \E[P^v_{\omega,\beta
}(T_{\parent v} = \infty)]$. Note that $f$ does not depend on $v$ by
the definition of the measure $\P$.
As mentioned in Section~\ref{sec:definition}, under assumptions
(XM)\vadjust{\goodbreak}
and (XR), there exists almost surely two infinite rays $v_0,\ldots
,v_n,\ldots$ and $w_0,\ldots,w_n,\ldots$ with $\liminf_{n\to\infty
}S(v_n)/n>0$ and $\limsup_{n\to\infty}S(w_n)/n<0$. By \eqref
{eq:conductance_ratio_est}, the former has finite resistance if $\beta
> -\beta_0$, and the latter if $\beta< -\beta_0$, whence $f(\beta) >
0$ for each $\beta\ne-\beta_0$. If $\beta= -\beta_0$, the Metropolis
algorithm is just a simple random walk on the $d$-regular tree and
therefore $f(-\beta_0) > 0$ as well. It follows that $f$ is positive
for every $\beta\in\R$. Furthermore, $f(\beta)$ is continuous because
it is the decreasing limit as $L\to\infty$ of $\E[P^v_{\omega,\beta
}(T_L < T_{\parent v})]$ and each of these quantities depends only on a
finite portion of the tree and is therefore continuous in $\beta$ by
assumption (H2). This immediately implies the lemma.
\end{pf}

The following important lemma controls quenched hitting probabilities
and will be used
in the evaluation of quenched escape probabilities.
%
%
\begin{lemma}
\label{lem:level_hitting_prob}
For each $K>0$, there exist $c,L_0>0$, $b>1$ depending on $K$,
such that for $|\beta|\leq K$ and $L>L_0$,
\[
\P\bigl(\exists v\in[1,L-1]\dvtx P_{\omega,\beta}^v(T_L
< T_{\parent
v}) < L^{-c} \bigr) < e^{-L^b}.
\]
\end{lemma}
\begin{pf}
Let $v\in[1,L-1]$, and let $u\in[L-|v|]_v$, such that $|u|=L$. By
\eqref
{eq:conductance_ratio_est}, we have for a fixed environment $\omega$,
for some $c>0$,
\[
P_{\omega,\beta}^v(T_u < T_{\parent v}) = \biggl(
\sum_{v\le w\le u} \frac{C(\parent v,v)} {C(\parent w,w)} \biggr)^{-1}
> c \biggl(\sum_{v\le w< u} e^{(\beta_0+\beta)(S(w)-S(v))}
\biggr)^{-1}.
\]
This gives for $|\beta|\le K$,
\begin{eqnarray*}
P_{\omega,\beta}^v(T_L < T_{\parent v})& \ge&\max
_{u\in[L-|v|]_v} P_{\omega,\beta}^v(T_u <
T_{\parent v})\\
 &>& \frac{c} L \max_{u\in
[L-|v|]_v} \min
_{w\in[v,u)} e^{-(K+\beta_0)|S(w)-S(v)|}.
\end{eqnarray*}
The statement now follows from Lemma~\ref{lem:brw_all_good}.
\end{pf}

For a vertex $v$, denote by $\ell(v)$ the depth of the first excursion
below $v$ after $T_v$, that is,
%
%
\begin{eqnarray}
\label{eq-starend} \ell(v) &=& \sup\bigl\{|V_n|-|v|\dvtx n\ge
T_v \mbox{ and }|V_k| > |v|\ \forall k\in\{T_v+1,
\ldots,n\} \bigr\}
\nonumber
\\[-8pt]
\\[-8pt]
\nonumber
&\in&\N\cup\{\infty\}.
\end{eqnarray}
Note that since the probability of jumping from $v$ to one of its
children does not involve $X(\parent v,v)$, the event $\ell(v) > 0$ is
independent from $X(\parent v,v)$ (conditioned on $T_v < \infty$).
%
%
\begin{lemma}
\label{lem:excursion_tail}
For each $K>0$ there exist $\alpha=\alpha(K),L_0=L_0(K) > 0$,
such that for $|\beta|\le K$ and
$L>L_0$,
\[
\P\bigl(P_{\omega,\beta}\bigl(L\le\ell(\root) < \infty\bigr)
> e^{-\alpha L}
\bigr) < e^{-\alpha L}.
\]
\end{lemma}

\begin{pf}
Fix $\ep\in(0,1)$. For a vertex $v$, define the variable $A(v)$ by
$A(v) = 1$ if for one of $v$'s
sisters $\widebar v$,
one has $P^{\widebar v}_{\omega,\beta}(T^*_{\widebar v} = \infty,
T_{\parent
v}=\infty) \ge\ep$ and
$A(v) = 0$ otherwise. In words, $A(v)=1$ if the
(quenched) probability of a walk, started at an appropriate sister
$\widebar v$
of $v$, to
escape to infinity through the subtree rooted at $\widebar v$
without visiting again $\widebar v$ is at least $\ep$.
By Lemma~\ref{lem:transience}, we
can choose $\ep$ such that $\E[A(v)] > 1/2$ for all $|\beta|\le
K$. By a result due to Grimmett and Kesten
(see \cite{Dembo2002}, Lemma~2.2, (2.1) for this version),
there exist then $\alpha,\gamma> 0$, such that
$\PP_\beta(G_L) \ge1-e^{-\alpha L}$ for large $L$, where
\[
G_L= \biggl\{\min_{v\in[L]} \sum
_{w\in[2,v]} A(w) \ge\gamma L \biggr\}.
\]
%
Now, let $\omega\in G_L$. We wish to bound $P_{\omega,\beta}(L\le
\ell
(\root) < \infty)$. For this, define $T^m$ for $m=2,\ldots,L-1$ to be
the first time after $T_L$ that the Markov chain $(V_n)_{n\ge0}$ hits
level $m$. If $T^m < \infty$ and $A(V_{T^m}) = 1$, then by assumptions
(XS) and (H1)--(H3), the probability that from $V_{T^m}$ the chain reaches
$\widebar{V}_{T^m}$ after two steps is bounded from below by $\delta
/\ep$ for some $\delta$ sufficiently small. It follows that
for $\omega\in G_L$,
\[
P_{\omega,\beta}\bigl(L\le\ell(\root) < \infty\bigr) \le(1-
\delta)^{\sum
_{w\in
[2,V_T]} A(w)} \le(1-\delta)^{\gamma L}.
\]
This yields the lemma (reducing the value of $\alpha$ if necessary).
\end{pf}

%
\begin{lemma}
\label{lem-alpha}
For each $K>0$ there exist $\alpha=\alpha(K),L_0=L_0(K) > 0$,
such that for $|\beta|\le K$ and
$L>L_0$,
\[
\PP_\beta\bigl(|V_{\tau_1}| \ge L\bigr) < e^{-\alpha L}.
\]
\end{lemma}
\begin{pf}
Define a sequence of random numbers $L_0,L_1,\ldots$ recursively as follows:
\begin{itemize}
\item$L_0 = 1$;
\item for $n\in\N$, let $v_n = V_{T_{L_n}}$. If $\ell(v_n) < \infty$,
then $L_{n+1} = L_n+\ell(v_n)+1$;
\item otherwise, set $L_m = \infty$ for $m>n$.
\end{itemize}
Let $N$ be the largest number $n$, such that $L_n < \infty$. Then by
construction, $|V_{\tau_1}| = L_N$. Furthermore, the differences
$(L_{n+1}-L_n)_{0\le n< N}$ are independent and identically distributed
as $\ell+1$ conditioned on $\ell<\infty$, and $N$ is geometrically
distributed with
success probability $\P(\ell= \infty) > 0$
[$\ell$ as in \eqref{eq-starend}].
The lemma then follows from Lemma~\ref{lem:excursion_tail}.
\end{pf}

%
\begin{lemma}
\label{lem:fresh}
Let $\mathcal G_n$ be the $\sigma$-field generated by $V_0,\ldots,V_n$
and let $T$ be a stopping time with respect to the filtration
$(\mathcal G_n)_{n\ge0}$, such that $V_T\ne V_k$ for all $k<T$.
Then for each $K>0$ there exists a constant $c=c(K)>0$,
such that for $|\beta|\le K$ and all $N\ge0$, we have
\[
\PP_\beta\Bigl(\max_{T \le j < T+N} |V_j| \ge cN
\big| \mathcal G_T \Bigr) > c.
\]
\end{lemma}
\begin{pf}
We follow the proof of \cite{Aidekon2008}, Theorem~1.5. Throughout the proof,
$c_0,c_1,\ldots$ will denote positive constants which are uniform in
$|\beta|\leq K$.

\emph{Step} 1. For a vertex $v\ne\root$, define $\pi_{\omega,\beta}(v)
= P^v_{\omega,\beta}(T_{\parent v} = \infty)$. Note that the random
variables $\pi_{\omega,\beta}(v)$, $v\ne\root$, are identically
distributed under $\P$ (but not independent). Let $\pi_{\omega,\beta}$
denote a random variable with this law. We wish to show that for some
constant $c_0$,
%
%
\begin{equation}
\label{eq:pi_omega_beta} \E[1/\pi_{\omega,\beta}] \le c_0.
\end{equation}
Denote by $v_1,\ldots,v_{d-1}$ the children of the vertex $v$. By
assumptions (XS) and (H1)--(H3), we have
%
%
\begin{equation}
\label{eq:one_step} P^v_{\omega,\beta}(V_1 =
v_i) \ge c_1,\qquad i=1,\ldots,d-1,
\end{equation}
which yields $\pi_{\omega,\beta}(v) \ge c_1 \sum_{i=1}^{d-1} \pi
_{\omega
,\beta}(v_i)$. Now, note that the variables $\pi_{\omega,\beta}(v_i)$,
$i=1,\ldots,i-1$ are independent under $\P$. The previous inequality
then yields for every $x\ge0$,
%
%
\begin{eqnarray}
\label{eq:pi_omega_beta_1} \P\bigl(\pi_{\omega,\beta}(v) \le
x\bigr) &\le&\P\Bigl(\max
_{i=1,\ldots,d-1} \pi_{\omega
,\beta}(v_i) \le
x/c_1\Bigr)
\nonumber
\\[-8pt]
\\[-8pt]
\nonumber
& =& \P\bigl(\pi_{\omega,\beta}(v) \le x/c_1
\bigr)^{d-1}.
\end{eqnarray}
Furthermore, by Lemma~\ref{lem:transience}, there exists a constant
$c_2$, such that
%
%
\begin{equation}
\label{eq:pi_omega_beta_2} \P\bigl(\pi_{\omega,\beta}(v) \le
2c_2\bigr) \le1/2.
\end{equation}
Together with \eqref{eq:pi_omega_beta_1}, this now easily implies
\eqref{eq:pi_omega_beta}.

\emph{Step} 2. For a vertex $v$, let $N_{v}$ denote the number of times
the vertex $v$ has been visited by the Metropolis algorithm
$(V_n)_{n\ge0}$. We wish to show that for each $k \ge0$,
%
%
\begin{equation}
\label{eq:Nv} \EE_\beta\biggl[\sum_{|v|=k}
N_{v} \biggr] \le c_3.
\end{equation}
Recall that $T^*_v$ denotes the first strictly positive hitting time of
the vertex $v$, such that $E^v_{\omega,\beta}[N_{v}] = 1/P^v_{\omega
,\beta}(T^*_v = \infty)$. By \eqref{eq:one_step}, we have
$P^v_{\omega
,\beta}(T^*_v = \infty)\ge c_1 \pi_{\omega,\beta}(v_1)$, such that
%
%
\begin{eqnarray}
\label{eq:Nv1} \EE_\beta[N_{v}] &=& \E\bigl[P^\root_{\omega,\beta}(T_v<
\infty)E^v_{\omega
,\beta}[N_{v}] \bigr] \le
c_1 \PP_\beta(T_v<\infty)\E\bigl[1/\pi
_{\omega
,\beta}(v_1)\bigr]
\nonumber
\\[-8pt]
\\[-8pt]
\nonumber
&\le& c_4
\PP_\beta(T_v<\infty),
\end{eqnarray}
by \eqref{eq:pi_omega_beta}. Furthermore, we have for every $k\ge0$,
%
%
\begin{eqnarray}
\label{eq:Tv} 1 &\ge&\sum_{|v|=k}\PP_\beta
(T_v<\infty, V_n \ge v_1\ \forall
n>T_v )
\nonumber
\\[-8pt]
\\[-8pt]
\nonumber
&\ge& c_1 \sum_{|v|=k}
\PP_\beta(T_v<\infty)\E\bigl[\pi_{\omega
,\beta}(v_1)
\bigr],
\end{eqnarray}
and by \eqref{eq:pi_omega_beta_2}, we have $\E[\pi_{\omega,\beta}(v_1)]
= \E[\pi_{\omega,\beta}] \ge c_2$. Equations~\eqref{eq:Tv} and
\eqref
{eq:Nv1} now yield \eqref{eq:Nv}.

\emph{Step} 3.
Recall the notation in the statement of the lemma, and let $L\in\N$.
Define the event $E_T = \{V_n \ge V_T\ \forall n>T\}$. A straightforward
extension of the proof of the last step allows us to prove that for
every constant $C>0$,
\[
\PP_\beta\biggl(E_T,\sum_{w\in[V_T,L]}
N_{w}>CL \Big| \mathcal G_T \biggr) \le\frac{1} {CL}
\EE_\beta\biggl[\sum_{w\in[V_T,L]}
N_{w}\Ind_{E_T} \Big| \mathcal G_T \biggr] <
c_3C^{-1}.
\]
Furthermore, by \eqref{eq:one_step}, we have $\PP_\beta(E_T |
\mathcal
G_T) \ge c_1 \E[\pi_{\omega,\beta}] \ge c_1c_2$. This now yields for
every constant $c>0$ and every $N\in\N$,
\[
\PP_\beta\Bigl(\max_{T \le j < T+N} |V_j| \ge cN
\big| \mathcal G_T \Bigr) \ge\PP_\beta\biggl(E_T,
\sum_{w\in[V_T,cN]} N_{w}\le N \Big| \mathcal
G_T \biggr)> c_1c_2 - c_3 c.
\]
Setting $c=c_1c_2/2c_3$ yields the proof.
\end{pf}

\def\Lb{\widebar L}

\begin{pf*}{Proof of Proposition~\ref{prop:tau_moments}}
The proof
follows an argument in the spirit of~\cite{Peres2008}, Proposition~3.
Let $K>0$.
Let $L=L(n)$ go to infinity with $n$ (we will later choose $L = n^{b}$
for some constant $b$).
We have, with $T_L$ denoting the hitting time
of level $L$, and with $\alpha>0$ as in the statement of Lemma~\ref
{lem-alpha},
%
%
\begin{eqnarray}
\label{eq-ofer1} \PP_\beta(\tau_1>n)& \le&
\PP_\beta\bigl(T_{|V_{\tau_1}|}>n, |V_{\tau
_1}| < L\bigr) +
\PP_\beta\bigl( |V_{\tau_1}|\ge L\bigr)
\nonumber
\\[-8pt]
\\[-8pt]
\nonumber
&\leq&\PP_\beta
(T_L>n)+e^{-\alpha L},
\end{eqnarray}
where in the last inequality we used the fact that $T_{L'} \le T_{L}$
for $L' \le L$.

Let $c$ be the constant from Lemma~\ref{lem:level_hitting_prob} and
set $\bar c=6 (c\vee1)$. Throughout the proof, all constants will be
uniform in $\beta$ for $|\beta| \le K$.
We write $\Lb= L^{\bar c}$.
Define a vertex $v$ to be \emph{fresh}
if it is visited by the random walk for the first time before
time~$\Lb$. We will upper bound $\PP_\beta (T_L>n)$ by showing that on
the one hand, there cannot be
too few fresh points that are well separated and on the other hand,
if there are many such fresh points, it is unlikely that $T_L$ is
large.

For a vertex $v$,
let $N_{v}$ denote the number of visits to $v$ by time $\Lb$, and
let $N_{v}^u$ denote the number of times the random walk visits $v$ before
time $\Lb$, and then,
at the next step that it moves, it visits the ancestor of $v$ before
time $\Lb$.
Clearly, for each $v$, $N_{v}^u\leq N_{v}$,
while,
using the
Markov property, there exist constants $\gamma,\gamma'>0$ such that
%
%
\begin{equation}
\label{eq-ofer3} \PP_\beta\bigl(\exists v\dvtx N_{v}\geq
\Lb^{1/2}, N_{v}^u\leq\gamma\Lb^{1/2}
\bigr) \leq\Lb e^{-\gamma' \Lb^{1/2}}.
\end{equation}
On the other hand,
the event
\[
\bigl\{N_{v}\geq\Lb^{1/2}, N_{v}^u>
\gamma\Lb^{1/2}, T_L\geq\Lb\bigr\}
\]
implies that the random walk visited
$v$ and then hit $\parent v$ at least $\gamma\Lb^{1/2}$ times.
By Lemma~\ref{lem:level_hitting_prob} and the Markov property, the
probability that there exists a fresh vertex $v$ satisfying
the last event
is bounded above by
\[
\Lb\bigl(1-L^{-c}\bigr)^{\gamma\Lb^{1/2}} \leq\Lb e^{-\gamma
L^{\bar c/2-c}} \leq
\Lb e^{- \gamma L^{2(c\vee1)}}.
\]
Combining this with
\eqref{eq-ofer3}, we conclude that for all large $L$,
%
%
\begin{equation}
\label{eq-ofer4} \PP_\beta\bigl(\exists v\dvtx N_{v}\geq
\Lb^{1/2}, T_L\geq\Lb\bigr) \leq e^{-L}.
\end{equation}

On the event $\{\forall v\dvtx N_{v}<\Lb^{1/2}\}$, there are at least
$\Lb
^{1/2}$ fresh points,
and therefore there are at least $\Lb^{1/4}$ fresh points that are
$\Lb^{1/4}$-separated. Let $c$ be the constant from Lemma~\ref
{lem:fresh}. At each arrival to such a fresh point, with
(annealed) probability at least $c $ the walk hits level $L$ before
time $L/c<\Lb^{1/4}$, for large $L$.
It follows that
%
%
\begin{equation}
\label{eq-ofer5} \PP_\beta\bigl(\forall v\dvtx N_{v}<
\Lb^{1/2}, T_L\geq\Lb\bigr) \leq C_0^{ \Lb^{1/4}}
\end{equation}
for some $C_0\in(0,1)$. By \eqref{eq-ofer4} and \eqref{eq-ofer5},
there exists now a constant $b>0$, such that with $L = n^b$, we have
$\PP_\beta(T_L > n) \le2e^{-L}$. Together with \eqref{eq-ofer1}, this
completes the proof of
the first statement of Proposition~\ref{prop:tau_moments}.

As for the second statement, let $\rho_{1},\ldots,\rho_{d-1}$ be the
children of the root. The law of the subtree of the branching random
walk tree rooted at $V_{\tau_1}$ is equal in law under $\PP_\beta$ to
the subtree rooted at $V_1$, conditioned on $V_1\in\{\rho_1,\ldots
,\rho
_{d-2}\}$ and on $T^*_\rho= \infty$. In particular,
\[
\PP_\beta(\tau_2-\tau_1\geq n)=
\PP_\beta\bigl(\tau_1\geq n | V_1\in\{
\rho_1,\ldots,\rho_{d-2}\}, T^*_\rho=\infty
\bigr) \le\PP_\beta(\tau_1\geq n)/c_{1}c_{2},
\]
where $c_1$ and $c_2$ are the constants from the proof of Lemma~\ref
{lem:fresh}. This finishes the proof.
\end{pf*}

\begin{pf*}{Proof of Theorem~\ref{th:einstein}}
Existence and finiteness of the limit $\sigma^2 = \lim_{n\to\infty}
S_n^2/n$ follows from the first part of Theorem~\ref{th:weak_er}. In
order to prove the remaining statements, we will use the regeneration
structure established in this section. Note that
the random vectors $\{(\tau_{i+1}-\tau_i, S_{\tau_{i+1}}-S_{\tau
_i})\}
_{i\geq1}$ are i.i.d. under
the law $\PP_\beta$ and independent from $(\tau_{1}, S_{\tau_{1}})$.
Furthermore, by (XS), $|S_{n}-S_{m}| \le g|n-m|$, where $g =
\operatorname{esssup}|X|$. Proposition~\ref{prop:tau_moments} and the
$\PP$-almost sure convergence of $S_n^2/n$ to $\sigma^2$ now yields
\[
\sigma^2 = \lim_{k\to\infty} \frac{S_{\tau_k}^2}{k}\cdot
\frac{k} {
\tau
_k} = \frac{\EE[(S_{\tau_2}-S_{\tau_1})^2]}{\EE[\tau_2-\tau_1]}
> 0,
\]
which proves the first part of the theorem.

As for the second part, by standard arguments (see, e.g., the proof
of Theorem~4.1 in \cite{Lyons1996}), $v_\beta= \lim_{n\to\infty}
S_n/n$ exists almost surely if $\EE_\beta[\tau_2-\tau_1]<\infty$, which
is the case for all $\beta$ by Proposition~\ref{prop:tau_moments}.
Furthermore, we have $v_\beta= \EE_\beta[S_{\tau_2}-S_{\tau_1}]
/\EE
_\beta[\tau_2-\tau_1]$, which implies in particular that $|v_\beta|
\le g$.

Now let $K_\beta= \inf\{k>0\dvtx\tau_k > \beta^{-2}\}$. By the optional
stopping theorem, we have
\[
\EE_\beta[S_{\tau_{K_\beta}}] = \EE_\beta[S_{\tau_1}] +
v_\beta\EE_\beta[\tau_{K_\beta}-\tau_1],
\]
such that by Proposition~\ref{prop:tau_moments} and assumption (XS),
for some constant $C>0$,
\[
\bigl|\EE_\beta[S_{\lfloor\beta^{-2}\rfloor}] - v_\beta\beta^{-2}\bigr|
\le C \EE_\beta\bigl[\tau_{K_\beta}-\beta^{-2}\bigr].
\]
Crude moment bounds using Proposition~\ref{prop:tau_moments}
yield that the right-hand side of the above equation is $o(\beta
^{-1})$, which yields $\lim_{\beta\to0} v_\beta/\beta= \lim
_{\beta\to
0} \EE_\beta[S_{\lfloor\beta^{-2}\rfloor}] = \sigma^2/2$ by
Theorem~\ref
{th:weak_er}. This finishes the proof.
\end{pf*}

\section*{Acknowledgments} We thank David Aldous for bringing the
problem of the analysis of the Metropolis algorithm on the tree
to our attention. We further thank an anonymous referee for
detailed comments improving the exposition of the article.


%
%



\printaddresses

\end{document}